\documentclass[12pt]{article}
\begin{document}
\begin{center}
{\large {\bf  IRREDUCIBLE REPRESENTATIONS OF 
U$_{pq}$[gl(2/2)]}} 
\vskip .8truecm
{\bf Nguyen Anh Ky}\\[2mm]
{\it Institute of Physics, P.O. Box 429, Bo Ho, 
Hanoi 10000}, Vietnam\\[1mm]
\end{center}
\vspace*{2mm}
\begin{abstract}     
The two--parametric quantum superalgebra 
$U_{p,q}[gl(2/2)]$ and its representations are 
considered. All finite--dimensional irreducible 
representations of this quantum superalgebra 
can be constructed and classified into typical 
and non--typical ones according to a proposition 
proved in the present paper.  This proposition is a 
nontrivial deformation from the one for the classical 
superalgebra $gl(2/2)$, unlike the case of 
one--parametric deformations.
\end{abstract}

\vspace*{4mm}

 
\underline{MSC--class.} : 81R50; 17A70.\\[1cm]
{\large {\bf  I. Introduction}}\\

  The quantum groups \cite{frt}--\cite{woro} were  introduced 
in 80's as a result of the study on quantum integrable systems 
and Yang--Baxter equations (YBE's)~\cite{collect}. It turns out 
that they are related to unrelated, at first sight, areas of both 
physics and mathematics and therefore,  have been intensively 
investigated in various aspects including their applications 
(see Refs.~\cite{frt}--\cite{ky1} and references therein). 
For applications of quantum groups, as in the non--deformed cases, 
we often need their explicit representations, in particular, the 
finite--dimensional ones which in many cases are connected 
with rational and trigonometric solutions of the quantum YBE's
\cite{frt}--\cite{2dim}. However, in spite of efforts and 
remarkable results in this direction the problem of investigating 
and constructing explicit 
representations of quantum groups, especially those for 
quantum superalgebras, is still far from being satisfactorily 
solved.  Even in the case of one--parametric quantum superalgebras, 
explicit representations are mainly known for quantum 
Lie superalgebras of lower ranks and of particular types like 
$U_{q}[osp(1/2)]$ and  $U_{q}[gl(1/n)]$ (Refs. \cite{ky1,cpt,pt}).  
So far, finite--dimensional representations of some bigger quantum 
superalgebras such as $U_{q}[osp(1/2n)]$ and $U_{q}[gl(m/n)]$  
with $m,n> 2$ have been considered but have not been explicitly 
constructed (see, for example, \cite{zhang, zhang2}). At the moment, 
detailed results in this aspect are known only for the cases with 
both $m,n\leq 2$ considered in \cite{ky1,ky2,ky3}, while for 
$U_q[gl(m/n)]$ with arbitrary $m$ and $n$ not all 
finite--dimensional representations but only a, although big, 
class of representations called essentially typical is known 
\cite{pnj}.\\

  As far as the multi--parametric deformations (first considered in 
\cite{man}) are concerned, this area is even less covered and 
results are much poorer.  Some kinds of two--parametric deformations 
have been considered by a number of authors from different points of 
view (see \cite{ky4}, \cite{jmp41} and references therein) 
but, to our knowledge, explicit representations are known and/or 
classified in a few lower rank cases such as 
$U_{p,q}[sl(2/1)]$ and $U_{p,q}[gl(2/1)]$ only 
\cite{ky4, zhang3}. The latter two--parametric quantum 
superalgebra $U_{p,q}[gl(2/1)]$ was consistently introduced 
and investigated in \cite{ky4} where all its finite--dimensional 
irreducible representations were explicitly constructed and 
classified at generic deformation parameters. This 
$U_{p,q}[gl(2/1)]$, however, is still a small quantum 
superalgebra which can be defined without the so--called 
extra--Serre defining relations \cite{extra1,extra2,extra3}
representing additional constraints on odd Chevalley generators 
in higher rank cases.  In order to include the extra--Serre 
relations on examination we introduced and considered a 
bigger two--parametric quantum superalgebra,  namely 
$U_{p,q}[gl(2/2)]$, and its representations \cite{jmp41,doson97}. 
Another our motivation for considering  this quantum superalgebra 
 is that already in the non--deformed case, 
the superalgebras $gl(n/n)$, especially, their subalgebras $sl(n/n)$ 
and $psl(n/n)$, have special properties (in comparison with other 
$gl(m/n)$, $m\neq n$) and, therefore, attract interest 
\cite{sigma1, sigma2,sigma3}. 
Additionally, structures of two--parameter deformations investigated 
in \cite{ky4, jmp41,doson97} and here are, of course, richer than those 
of one--parameter deformations. Every deformation parameter can be 
independently chosen to take a separate generic value (including zero) 
or to be a root of unity.\\
     
  Combining the advantages of  the previously developed 
methods \cite{ky2, ky3,ky4} for $U_{q}[gl(2/2)]$ and 
$U_{p,q}[gl(1/2)]$ we  described in \cite{jmp41} how to 
construct finite--dimensional representations of  the 
two--parametric quantum Lie superalgebra $U_{p,q}[gl(2/2)]$.
In this paper  we consider when these representations constructed 
are irreducible.  It turns out that they can be classified again 
into typical and nontypical representations which, even at generic 
deformation parameters, however, are nontrivial deformations from 
the classical analogues \cite{kkp}, unlike many cases of 
one-parametric deformations.\\[7mm] 
{\large {\bf  II. The quantum superalgebra $U_{p,q}[gl(2/2)]$}}\\ 

   The quantum superalgebra $U_{p,q}\equiv U_{p,q}[gl(2/2)]$
as a two--parametric deformation of the universal enveloping 
algebra $U[gl(2/2)]$ of the Lie superalgebra $gl(2/2)$ can be 
completely generated by the operators 
$L_k$, $E_{12}$, $E_{23}$, $E_{34}$, 
$E_{21}$, $E_{32}$, $E_{43}$ and $E_{ii}$ 
($1\leq i \leq 4$)
called again Cartan--Chevalley generators subjects to the 
following (defining) relations \cite{jmp41,doson97}: 
\begin{tabbing}
\=123456$[E_{12},E_{34}]$\=$[E_{ii},E_{jj}]$12345678912
\= =1\= 0,1234 \=$[E_{ii},E_{j,j+1}]$\= =
\=$(\delta_{ij}-\delta_{i,j+1})1234~$\=\kill

{}~~~~a) {\it the super--commutation relations} 
($1\leq i,i+1,j,j+1\leq
4$):\\[2mm]
\>\>$[E_{ii},E_{jj}]$\> = \>0,\>\>\>\>(1a) \\[1mm] 
\>\>$[E_{ii},E_{j,j+1}]$\>=\>$(\delta_{ij}-\delta_{i,j+1})E_{j,j+1},$ 
\>\>\>\>(1b)\\[1mm]
\>\>$[E_{ii},E_{j+1,j}]$\>=\>$(\delta_{i,j+1}-\delta_{ij})E_{j+1,j}$, 
\>\>\>\>(1c)\\[1mm]
     \>\>[even generator, $L_{k}$]\>=\>0,~ $k=1,2,3$,~
\>\>\>\>(1d)\\[1mm]
\>\>$[E_{i,i+1},E_{j+1,j}\}$\>=\>$\delta_{ij} 
\left({q\over p}\right)^{L_{i}-H_{i}(1+\delta_{i2})/2} 
[H_{i}]$,
\>\>\>\>(1e)
%\\[4mm]
\end{tabbing}
\begin{tabbing}
\=123456781$[E_{12},E_{34}]$=\=$[E_{ii},E_{jj}]$1234\= =1\= 0,1234 
\=$[E_{ii},E_{j,j+1}]$\= =
\=$(\delta_{ij}-\delta_{i,j+1})E_{j,j+1}$12~~\=\kill 
%\\

{}~~~~b) {\it the Serre--relations}:\\[2mm]
 \>\>~~$[E_{12},E_{34}]$\>=\>$[E_{21},E_{43}]$\>~~~~~~~=~~0,
\>\>\>(2a)\\[1mm]
\>\>~~~~~~~$E_{23}^{2}$\>=\>~~~~~$E_{32}^{2}$\>~~~~~~~=~~0, 
\>\>\>(2b)\\[1mm]
     
\>~~~~~~~~~~
$[E_{12},E_{13}]_{p}$~~=\> ~~$[E_{21},E_{31}]_{q}$\>=~~ 
\>$[E_{24},E_{34}]_{q}$\>~~~~~~~=~~$[E_{42},E_{43}]_{p}=~~0$,
                           \>\>\>(2c)
\end{tabbing}
and\\

c) {\it the extra--Serre relations}:\\ 
$$\{E_{13},E_{24}\}=0,\eqno(3a)$$
$$\{E_{31},E_{42}\}=0,\eqno(3b)$$
where $H_{i}\equiv (E_{ii}-{d_{i+1}\over d_{i}}E_{i+1,i+1})$, 
$d_{1}=d_{2}=-d_{3}=-d_{4}=1$, 
$L_{1}\equiv L_{l}, L_{2}\equiv 0, L_{3}\equiv L_{r}$ (with $L_l$ 
and $L_r$ explained later),    
$[x]\equiv (q^x-p^{-x})/(q-p^{-1})$ is a so-called $pq$--deformation 
of $x$ being a number or an operator and, finally, [ , \} is a notation 
for the supercommutators.  Here, the operators
$$E_{13}~:=~[E_{12},E_{23}]_{q^{-1}},\eqno(4a)$$ 
$$E_{24}~:=~[E_{23},E_{34}]_{p^{-1}},\eqno(4b)$$ 
$$~~~E_{31}~:=~-[E_{21},E_{32}]_{p^{-1}},\eqno(4c)$$ 
$$~~~E_{42}~:=~-[E_{32},E_{43}]_{q^{-1}}\eqno(4d)$$ 
and the operators composed in the following way 
\begin{tabbing}
\=123456789123456\=$E_{41}$~\=:=1
\=$[E_{21},[E_{32},E_{43}]_{q^{-2}}]_{q^{-2}}$
\=$\equiv [E_{21},E_{42}]_{q^{-1}}$1234567891234~~~~\=\kill 
\>\>$E_{14}$~\>:=\>$[E_{12},[E_{23},E_{34}]_{p^{-1}}]_{q^{-1}}$ 
\>$\equiv ~[E_{12},E_{24}]_{q^{-1}}$,\>(5a)\\[1mm] 
\>\>$E_{41}$\>:=\>$[E_{21},[E_{32},E_{43}]_{q^{-1}}]_{p^{-1}}$ 
\>$\equiv ~-[E_{21},E_{42}]_{p^{-1}}$\>(5b)
\end{tabbing}
are defined as new generators, where $[A,B]_r = AB-rBA$. 
%is a so-called $r$--deformed commutator. 
These generators, like $E_{23}$ and 
$E_{32}$, are all odd and have vanishing squares. 
The generators $E_{ij}$, $1\leq i,j\leq 4$, are
two--parametric deformation analogues ($pq$--analogues) of the 
Weyl generators $e_{ij}$ of the superalgebra 
$gl(2/2)$ whose universal enveloping algebra $U[gl(2/2)]$ is a 
classical limit of $U_{p,q}[gl(2/2)]$ when $p,q\rightarrow 1$. 
The so--called maximal--spin operator $L_l$ (or $L_r$)  
is a constant within a finite--dimensional irreducible module  
({\it  fidirmod}) of  a  $U_{p,q}[gl(2)]$ (defined below)   
and are different for different $U_{p,q}[gl(2)]$--fidirmods. 
Therefore, commutators between these operators with the odd 
generators intertwining $U_{p,q}[gl(2)]$--fidirmods take 
concrete forms on concrete basis vectors. Other commutation 
relations between $E_{ij}$ follow from the relations (1)--(3) 
and the definitions (4) and (5).\\[7mm]
{\large {\bf  III.  Representations of  $U_{p,q}[gl(2/2)]$}}\\
     
 The subalgebra $U_{p,q}[gl(2/2)_{0}]~ (\subset
U_{p,q}[gl(2/2)]_{0}\subset U_{p,q}[gl(2/2)])$ is even and 
isomorphic to $U_{p,q}[gl(2)\oplus gl(2)]\equiv U_{p,q}[gl(2)] 
\oplus U_{p,q}[gl(2)]$ which can be completely generated by 
$L_1$, $L_3$, $E_{12}$, $E_{34}$, $E_{21}$, $E_{43}$ and 
$E_{ii}$, $1\leq i\leq 4$,
$$U_{q}[gl(2/2)_{0}]~=~ 
\mbox{lin.env.} \{L_1, L_3, E_{ij}\|~ 
i,j=1,2~~
\mbox{and}~~i,j=3,4\}.\eqno(6)$$ 
In order to distinguish two components $U_{p,q}[gl(2)]$ 
of $U_{p,q}[gl(2/2)_{0}]$ we set 
$$\mbox{left}~~U_{p,q}[gl(2)]\equiv U_{p,q}[gl(2)_{l}]:=
\mbox{lin.env.}\{L_1, E_{ij}\|~i,j=1,2\},\eqno(7)$$ 
$$\mbox{right}~U_{p,q}[gl(2)]\equiv U_{p,q}[gl(2)_{r}]:= 
\mbox{lin.env.}\{L_3, E_{ij}\|~ i,j=3,4\}, 
\eqno(8)$$
that is
$$U_{p,q}[gl(2/2)_{0}]~=~U_{p,q}[gl(2)_{l}\oplus gl(2)_{r}]. 
\eqno(9)$$
     
%%%%%%
   We see that every of
the odd spaces $A_{+}$ and $A_{-}$ spanned on the positive 
and negative odd roots (generators) $E_{ij}$ and $E_{ji}$, 
$1\leq i\leq 2 <j\leq 4$, respectively
$$A_{+}= 
%{\normalsize 
\mbox{lin.env.}\{E_{14},E_{13}, 
E_{24}, E_{23}\},\eqno(10)$$
$$A_{-}= 
\mbox{lin.env.}\{E_{41},E_{31}, E_{42},E_{32} \}, 
\eqno(11)$$
is a representation space of the even subalgebra 
$U_{p,q}[gl(2/2)_{0}]$
which, as seen from (1)--(2), is a stability subalgebra 
of $U_{p,q}[gl(2/2)]$.
Therefore, we can construct representations of 
$U_{p,q}[gl(2/2)]$ induced from some (finite--dimensional 
irreducible, for example) representations of 
$U_{p,q}[gl(2/2)_{0}]$ which are realized
in some representation spaces (modules) $V^{p,q}_{0}$ 
being tensor products of
$U_{p,q}[gl(2)_l]$--modules $V^{p,q}_{0,l}$ and 
$U_{p,q}[gl(2)_r]$--modules $V^{p,q}_{0,r}$ 
$$V_{0}^{p,q}(\Lambda)=V_{0,l}^{p,q}(\Lambda_{l})\otimes 
V_{0,r}^{p,q}(\Lambda_{r}),\eqno(12)$$
where  $\Lambda$'s are some signatures (such as highest 
weights, respectively) characterizing the modules (highest 
weight modules, respectively). Here $\Lambda_l$ and 
$\Lambda_r$ are referred to as the left and the right 
components of $\Lambda$, respectively,
$$\Lambda=[\Lambda_{l},\Lambda_{r}].\eqno(13)$$ 
     
  If we demand
$$E_{23}V_{0}^{p,q}(\Lambda)=0 \eqno(14)$$ hence 
$$U_{p,q}(A_{+})V_{0}^{p,q}=0,\eqno(15)$$ we 
turn the $U_{p,q}[gl(2/2)_{0}]$--module 
$V^{p,q}_{0}$ into a $U_{p,q}(B)$--module where 
$$B=A_{+}\oplus gl(2)\oplus gl(2).\eqno(16)$$
The $U_{p,q}[gl(2/2)]$--module $W^{p,q}$ induced from 
the $U_{p,q}[gl(2/2)_{0}]$--module $V^{p,q}_{0}$ is the 
factor--space
$$W^{p,q}=W^{p,q}(\Lambda)=[U_{p,q}\otimes V_{0}^{p,q} 
(\Lambda)]/I^{p,q}(\Lambda)\eqno(17)$$
which, of course, depends on $\Lambda$, where 
$$U_{p,q}\equiv U_{p,q}[gl(2/2)], \eqno(18)$$ 
while $I^{p,q}$ is the subspace
$$I^{p,q}=
\mbox{lin.env.}
\{ub\otimes v-u\otimes bv\| u\in U_{p,q},
b\in U_{p,q}(B)\subset U_{p,q}, v\in V_{0}^{p,q}\}. 
\eqno(19)$$

   Using the commutation relations (1)--(3) and the
definitions (4) and (5) we can prove the an 
analogue of the Poincar\'e--Birkhoff--Witt  theorem. 
Consequently, a basis of 
$W^{p,q}$ can be constituted by taking 
all the vectors 
of the form $$\left |\theta_{1}, \theta_{2}, \theta_{3}, 
\theta_{4}; (\lambda)\right > :=
(E_{41})^{\theta_{1}}
(E_{31})^{\theta_{2}}(E_{42})^{\theta_{3}}(E_{32})^ 
{\theta_{4}}\otimes (\lambda),
~~ \theta_{i}=0,1, \eqno(20)$$ where $(\lambda)$ 
is a (Gel'fand--Zetlin, for example) basis of $V_{0}^{p,q}\equiv 
V_{0}^{p,q}(\Lambda)$. This basis of 
$W^{p,q}$ called the induced $U_{p,q}[gl(2/2)]$--basis
(or simply, the induced basis), however, is not convenient 
for investigating the module structure of $W^{p,q}$.  It 
was the reason the so-called reduced basis was introduced 
\cite{jmp41}. It is obvious that if the module $V_{0}^{p,q}$ 
is finite--dimensional so is the module $W^{p,q}$. 
In this case  $W^{p,q}$ can be characterized by a signature 
$[m]$ and is decomposed into a direct sum of (sixteen, at most) 
$U_{p,q}[gl(2/2)_{0}]$--fidirmod's  $V_{k}^{p,q}$ 
of  signatures $[m]_k$ :
$$W^{p,q}([m])=\bigoplus_{k=0}^{15}V_{k}^{p,q}([m]_{k}). 
\eqno(21)$$ 
Thus, the reduced basis of $W^{p,q}$ is a union of the bases 
of all $V_{k}^{p,q}$'s  which can be presented by the 
quasi--Gel'fand--Zetlin patterns \cite{jmp41}, 
corresponding to the branching rule 
$U_{p,q}[gl(2/2)]\supset U_{p,q}[gl(2/2)_{0}]\supset 
U_{p,q}[gl(1)\otimes gl(1)]$,
$$
\left[
\begin{array}{lccc}
     
                        m_{13}& m_{23}&  m_{33}& m_{43} \\
     
                        m_{12}& m_{22}&  m_{32}&  m_{42}\\
m_{11}&   0   &  m_{31}&    0
\end{array}
\right]_{k}
\equiv
(m)_{k}, ~~ 0\leq k\leq 15,
\eqno(22)$$
where $m_{ij}$ are complex numbers such that 
$m_{i2}-m_{i1}\in 
{\bf Z}^{+}$, $m_{i1}-m_{i+1,2}\in {\bf Z}^{+}$, 
$m_{i3}-m_{i+1,3}\in {\bf Z}^{+}$, $i=1,3.$
The  second row 
$[m_{12}, m_{22},m_{32},m_{42}]$ in (22)
is fixed for a given $k$, as for $k=0$ it takes the 
value of the first row $[m_{13}, m_{23},m_{33},m_{43}]$  
which is fixed for all $k=0,1,...15$.
Now, a  signature $[m]_k$ of a $V_k^{p,q}$ is identified with 
a second row, 
$$[m]_k\equiv [m_{12}, m_{22},m_{32},m_{42}],$$ 
while the signature $[m]$  single in the whole $W^{p.q}$  
(i.e., the same for all $V_k^{p.q}$'s) is indentified with 
the first row, $$[m]\equiv [m_{13}, m_{23},m_{33},m_{43}].$$ 
The actions of the generators 
$E_{ij}$ on the basis (22) are given in \cite{jmp41} or
can be calculated by using the method explained there. 
The basis vector (22) with $m_{11}=m_{12}$ and  $m_{31}=m_{32}$
$$ (M)_{k}=
\left[
\begin{array}{lccc}
     
                        m_{13}& m_{23}&  m_{33}& m_{43} \\
     
                        m_{12}& m_{22}&  m_{32}&  m_{42}\\
m_{12}&   0   &  m_{32}&    0
\end{array}
\right]_{k},
\eqno(23)$$
annihilated by $E_{12}$ and $E_{34}$ is, by definition, the 
highest weight vector of  the submodule $V_k^{p,q}([m]_k)$. 
For $k=0$ the highest weight vector of  the submodule 
$V_0^{p,q}([m])$
$$ (M)_0\equiv (M)=
\left[
\begin{array}{lccc}
     
                        m_{13}& m_{23}&  m_{33}& m_{43} \\
     
                        m_{13}& m_{23}&  m_{33}&  m_{43}\\
m_{13}&   0   & m_{33}& 0 \end{array} \right].  \eqno(24)$$ 
is, in addition, also annihilated by the odd genrator $E_{23}$ 
and, therefore,  simultaneously represents the highest weight 
vector of both $V_0^{p,q}([m])$ and $W^{p,q}([m])$. A monomial 
of the form 
$$
\left |\theta_{1}, \theta_{2}, \theta_{3}, \theta_{4}\right > :=
(E_{41})^{\theta_{1}}(E_{31})^{\theta_{2}}(E_{42})^{\theta_{3}}
(E_{32})^{\theta_{4}}, ~~  \theta_i =0,~1\eqno(25)$$ 
would shift a subspace $V_k^{p,q}$ to another subspace $V_l^{p,q}$ 
with $l>k$. So here we would call the former a higher (weight) 
subspace with respect to  the latter called a lower (weight) 
subspace.\\[7mm]
{\bf Proposition}: {\it The induced module $W^{p,q}[m]$ constructed 
is irreducible if and only if 
\begin{eqnarray*}
&&[h^0_2][h^0_1+h^0_2+1]\left\{-{q\over p}[h^0_2-1][h^0_3+1] + 
[h^0_2][h^0_3]\right\}\times \\
&&\left\{-q^{-h^0_2+1}p^{-h^0_3-1}[h^0_1+1]-q^{h^0_1}
\left( {q\over p}\right)
^{-h^0_2+1}[h^0_2-1][h^0_3+1] + q^{h^0_1}\left( {q\over p}\right)
^{-h^0_2}[h^0_2][h^0_3] \right.\\
&&\left. +{q\over p}\left( -q^{-h^0_2} + q^{-h^0_2-2}\right)
[h^0_3]\left( q^{h^0_1+1}  + {q^2\over p^2}[h^0_1]\right)\right\}\neq 0.
~~~~~~~~~~~~~~~~~~~~~~~~~~~~~~~(26) 
\end{eqnarray*} 
where  ~ 
$h^0_1=m_{13}-m_{23}, ~ h^0_2=m_{23}+m_{33}, ~h^0_3=m_{33}-m_{43}.$}\\

  The irreducible module $W^{p,q}$ constructed with keeping the 
condition (26) valid is called typical, otherwise, we say it is 
an indecomposable module.  In the latter case, however, there always 
exists a maximal invariant submodule $I_h^{p,q}$  
(of class $h$, $h=1,2,...$)  of  $W^{p,q}$ and the compliment to 
$I_h^{p,q}$ 
subspace of $W^{p,q}$ is not invariant under $U_{p,q}[gl(2/2)]$ 
transformations. The representation carried in the  factor module 
$W^{p,q}/I_h^{p,q}$ is irreducible and called a nontypical 
representation of  $U_{p,q}[gl(2/2)]$.  It can be shown that these 
typical and nontypical representations contain all classes of 
finite--dimensional irreducible representations of 
$U_{p,q}[gl(2/2)]$. \\ 

  As  every subspace  $V_{k}$, k=0,1,..., 15, is close and already 
irreducible under the even subalgebra $U_{p,q}[gl(2/2)_{0}]$, to 
see if  $W^{p,q}$ is an irreducible module of $U_{p,q}$ it remains 
to consider  the action of its odd generators only.  By construction 
(see Eqs. (17)--(21)) the module $W^{p,q}$ is at least indecomposable 
since any its subspace $V_{k}$, $1\leq k\leq 15$, including the 
lowest one $V_{15}$, can be always reached from higher subspaces 
$V_l$,  $0\leq l<k$, including the highest one $V_0$, acted by the 
monomials 
$\left |\theta_{1}, \theta_{2}, \theta_{3}, \theta_{4}\right >$ 
given in (25). Contrarily, the monomials 
$$\left <\theta_{1}, \theta_{2}, \theta_{3}, \theta_{4}\right |:=
(E_{14})^{\theta_{1}}(E_{13})^{\theta_{2}}(E_{24})^{\theta_{3}}
(E_{23})^{\theta_{4}}\eqno(27)$$ 
send us to the opposite direction: from lower subspaces to 
higher ones. Thus, the module  $W^{p,q}$ is irreducible if and 
only if $V_0$ is reachable from the lowest subspace $V_{15}$ 
under the action of the operators (27).  The most optimal way 
to see that is to act on a vector of the subspace $V_{15}$ by 
the monomial $E_{14}E_{13}E_{24}E_{23}$, i.e.,  
the monomial (27)  with all $\theta_i$'s  = 1 but not less 
(an action of a shorter monomial on $V_{15}$ should not reach 
$V_0$).  Since $V_{15}$ is an irreducible module of 
 $U_{p,q}[gl(2/2)_{0}]$, 
it is simplest but enough to consider when the highest 
weight vector  $E_{41}E_{31}E_{42}E_{32} (M)$ of 
$V_{15}$ under the action 
of  $E_{14}E_{13}E_{24}E_{23}$ reaches (or we can say, 
returns to) 
$V_{0}$. In other words, the module  $W^{p,q}$ is irreducible 
if and only if the condition 
$$E_{23}E_{24}E_{13}E_{14}E_{41}E_{31}E_{42}E_{32} 
(M)\neq 0 
\eqno(28)$$
holds. This condition in turn can be proved (for  $p,q\neq 0$) 
to be 
equivalent to the condition
\begin{eqnarray*}
& &
[H_2][H_1+H_2+1]\left\{-{q\over p}[H_2-1][H_3+1] + [H_2][H_3]
\right\}\times \\
&&
\left\{-q^{-H_2+1}p^{-H_3-1}[H_1+1]-q^{H_1}\left( {q\over p}\right)
^{-H_2+1}[H_2-1][H_3+1] + q^{H_1}\left( {q\over p}\right)
^{-H_2}[H_2][H_3] \right.\\
&&\left. +{q\over p}\left( -q^{-H_2} + q^{-H_2-2}\right)
[H_3]\left( q^{H_1+1}  + {q^2\over p^2}[H_1]\right)\right\}(M)\neq 0.
~~~~~~~~~~~~~~~~~~~~~~~~~(29) 
\end{eqnarray*} 
which is nothing but (26) with $h^0_i$ being eigenvalues of $H_i$ 
on the highest weight vector $(M)$. The proposition is, thus,  
proved.\\[7mm]
{\large \bf  IV. Conclusion}\\

  The two-parametric quantum superalgebra $U_{p,q}[gl(2/2)]$ 
was introduced in \cite{jmp41, doson97}. Its representations 
constructed by the method described in \cite{jmp41} are either 
irreducible (when the condition (26) is kept) or indecomposable 
(when the condition (26) is violated). The irreducible 
representations in the former case are called typical. In the 
case of indecomposable representations, however, irreducible 
representations can be always extracted. One such irreducible 
representation called nontypical is simply a 
factor--representation in a factor--subspace  of the original 
indecomposable module factorized by its maximal invariant 
subspace.  All the typical and nontypical representations 
constructed in such a way contain all classes of 
finite--dimensional irreducible representations of 
$U_{p,q}[gl(2/2)]$. For conclusion, let us emphasize 
that  the condition (26) and the representations become more 
interesting at roots of unity but they, even at generic 
deformation parameters, are nontrivial deformations from 
the classical analogues \cite{kkp} in the sense that the former 
cannot be found from the latter by replacing in appropriate 
places the ordinary brackets with the quantum deformation 
ones, unlike many one--parametric cases. \\[7mm]
{\bf Acknowledgements}\\

    This work was supported in part by the National Research
Programme for Natural Sciences of Vietnam under grant  
number KT -- 04.1.2.

\end{document}